# BINOMIAL PREDICTION USING THE FREQUENT OUTCOME APPROACH

B. O'NEILL,* *University of New South Wales*

## Abstract

Within the context of the binomial model, we analyse sequences of values that are almost-uniform and we discuss a prediction method called the frequent outcome approach, in which the outcome that has occurred the most in the observed trials is the most likely to occur again. Using this prediction method we derive probability statements for the prior probability of correct prediction, conditional on the underlying parameter value in the binomial model. We show that this prediction method converges to a level of accuracy that is equivalent to ideal prediction based on knowledge of the model parameter.

Secondary 62C10

## 1. Introduction

It is well known that exchangeable sequences of indicator values lead us to a binomial model. This model is used in gambling processes such as coin flipping, where there are only two possible outcomes in each trial. In such cases, outcomes are generated by repeated trials of some physical process under conditions which are identical in all relevant respects from trial to trial, leading to a judgment of infinite exchangeability. The decision problem confronting the gambler is to predict the outcomes of future trials of the process under consideration. In short, in simple gambling processes with two possible outcomes in each trial, and unchanging conditions, we are concerned with making predictions from a binomial model.

In predicting the outcomes of any exchangeable sequence of values, it can be shown that the values cannot be negatively correlated and will only be positively correlated when the long-run empirical distribution is treated as a random variable (see Section 2 for details). This is usually the case in Bayesian analysis, so that the outcomes of an exchangeable sequence are found to be positively correlated.

In fact, this positive correlation between repeated outcomes leads us naturally to an optimal prediction method called the *frequent outcome approach*, whereby the most frequent outcome is taken as the most likely outcome to occur, on the basis that the observed data indicate the likely direction of any underlying biases in the process. This approach is discussed by Wilson (1965) and is given formal mathematical treatment using Bayesian modelling in O'Neill and Puza (2005). The latter showed that the frequent outcome approach maximises the posterior probability of correct prediction for any exchangeable sequence of values, where the long-run proportions of outcomes are also exchangeable *a priori*. In the case of gambling processes, this leads to an optimal prediction method which holds that the most frequently occurring outcome is the most likely to occur again.

While this approach can be shown to maximise the posterior probability of correct prediction (under the conditions stated), it is not generally sufficient to solve the kinds of decision problems that arise in gambling situations. In these problems, in addition to the determination of the most likely outcome, the decision-maker also has the choice to gamble or not to gamble, which is affected by the *consequences* of correct and incorrect prediction. In this case, knowing that the posterior probability of correct prediction is maximised for a particular outcome is not enough – the decision-maker must be able to *evaluate* this probability (or at least approximate it) in order to know whether to gamble or not.

This is a similar problem to the one considered by Ethier (1982), who developed a hypothesis test to determine whether the bias in a gambling process is sufficient to obtain a bet having positive expected return using the frequent outcome approach. This test for 'favourable' outcomes is very valuable and allows the gambler to decide when to bet. However, it is also valuable to consider how long we would expect it to take for evidence of favourable bias to emerge, if such bias exists. To do this we need to look at the probability of correct prediction, conditional on the probability parameter in the binomial model.

In this paper we analyse the prediction problem that arises from an exchangeable series of values with two possible outcomes; that is, we analyse prediction under a simple binomial model. More specifically, we look at certain sequences that arise under the binomial model when a particular form of prior distribution is used to model the parameter. We look at the frequent outcome approach and assess the power of this approach conditional on given values of the underlying parameters in the model.

## 2. Exchangeability and the binomial model

Suppose that we have a sequence of indicator values $x := (x_1, x_2, x_3, \ldots)$ (i.e. each value is either zero or one). We define the associated *count values* $n := (n_1, n_2, n_3, \ldots)$ by $n_k := \sum_{i=1}^{k} x_i$ and we define the *long-run proportion* by $\theta := \lim_{k \to \infty} n_k / k$. We also note that the long-run proportion has range $0 \le \theta \le 1$.

If the sequence $x$ is exchangeable then it follows from the representation theorem of de Finetti (1980) that the elements of $x|\theta$ are independent Bernoulli random variables with occurrence probability $\theta$ (see also Fortini *et al.* (2000)). This in turn implies that

$$p(n_k \mid \theta) = \text{Bin}(n_k \mid k, \theta) = \binom{k}{n_k} \theta^{n_k} (1 - \theta)^{k - n_k}.$$

Thus, the judgment of exchangeability leads us to the standard binomial model where the unknown parameter of interest is the long-run proportion $\theta$. Under this model the count values form a sufficient statistic for inference about the long-run proportions.

O'Neill (2009) showed that the elements of an exchangeable sequence cannot have negative correlation, and will only have positive correlation if the long-run empirical distribution of the sequence is treated as a random variable (the first of these results can also be found in Kingman (1978)). In particular, if $x$ is exchangeable, then for all $i \ne j$ we have the following covariance:

$$\text{cov}(x_i, x_j) = \text{var}(\mu(\theta)) \ge 0, \quad \text{where } \mu(\theta) := \text{E}[x_i \mid \theta].$$

If $\theta$ is treated as a random variable (as in the Bayesian approach), and is sufficiently random to ensure that $\mu(\theta)$ has positive variance, then the elements of the sequence $x$ will be positively correlated. This positive correlation leads to a useful prediction method, which arises in problems displaying a form of informational symmetry.

### 3. The prediction problem

The prediction problem in the binomial model consists of specifying a method for predicting the next value of the observable sequence given previous observations of values in the sequence. To do this, suppose that we observe $x_1, \ldots, x_k$ for some $k \geq 0$ and we want to make a prediction for the next value, $x_{k+1}$.

Suppose that we let $\boldsymbol{y} := (y_1, y_2, y_3, \ldots)$ be our predicted values for the respective values in the sequence $\boldsymbol{x}$. Since the values in the sequence are exchangeable, we can make each prediction $y_{k+1}$ solely based on the sufficient statistic $(k, n_k)$ without any loss in predictive power. Thus, without loss of generality, we will have

$$y_{k+1} \sim \mathrm{Bin}(1, \phi_{k,n_k}), \quad \text{where } 0 \leq \phi_{k,n_k} \leq 1 \text{ for all } 0 \leq n_k \leq k.$$

This means that for any observed value $0 \leq n_k \leq k$ we will set our prediction as

$$y_{k+1} := \begin{cases} 0, & \text{with probability } 1 - \phi_{k,n_k}, \\ 1, & \text{with probability } \phi_{k,n_k}. \end{cases}$$

This general prediction method is fully characterized by the *prediction array* which is an array of numbers that looks like this:

$$\boldsymbol{\phi} = \begin{pmatrix} \phi_{0,0} & & & & \\ \phi_{1,0} & \phi_{1,1} & & & \\ \phi_{2,0} & \phi_{2,1} & \phi_{2,2} & & \\ \phi_{3,0} & \phi_{3,1} & \phi_{3,2} & \phi_{3,3} & \\ \vdots & \vdots & \vdots & \vdots & \ddots \end{pmatrix}.$$

Each of the numbers in the array gives a probability value for prediction; the rows correspond to different values of $k$ and the values in the rows correspond to different values of $n_k$. If the prediction probabilities in the array are all zero or one then this yields a deterministic prediction method; otherwise it yields a prediction method that can give nondetermined values in some instances.

Since the prediction array fully characterizes the prediction method, the prediction problem consists of specifying some prediction array which will be used to make our predictions. (If we only need to make a prediction for certain specified values of $k$ and $n_k$ then only some of the prediction array will need to be specified.)

Having specified this, the probability of correct prediction, conditional on both $n_k$ and $\theta$, is given by

$$\begin{aligned} \pi_{k,n_k}(\theta) &:= \mathrm{P}(y_{k+1} = x_{k+1} \mid n_k, \theta) \\ &= \mathrm{P}(y_{k+1} = x_{k+1} = 0 \mid n_k, \theta) + \mathrm{P}(y_{k+1} = x_{k+1} = 1 \mid n_k, \theta) \\ &= (1-\theta)(1 - \phi_{k,n_k}) + \theta \phi_{k,n_k} \\ &= 1 - \phi_{k,n_k} + (1 - 2\phi_{k,n_k})\theta. \end{aligned}$$

Of course, this equation is of little practical use on its own, given that we will not have knowledge of the parameter $\theta$. To obtain a sensible prediction method we will need to look at the posterior probability of correct prediction, *without* conditioning on the parameter in the model. (This is the situation faced by a gambler who is betting on the outcomes of an exchangeable sequence of values.)

### 4. The Bayesian approach to prediction

Under the Bayesian approach we specify a prior distribution for the parameter $\theta$ and we work conditionally on the observed value $n_k$. The posterior probability of correct prediction is given by

$$
\begin{aligned}
\mathrm{P}(y_{k+1} = x_{k+1} \mid n_k) &= \mathrm{E}[\pi_{k,n_k}(\theta) \mid n_k] \\
&= \mathrm{E}[1 - \phi_{k,n_k} + (1 - 2\phi_{k,n_k})\theta \mid n_k] \\
&= 1 - \phi_{k,n_k} + (1 - 2\phi_{k,n_k})\,\mathrm{E}[\theta \mid n_k] \\
&= \pi_{k,n_k}(\mathrm{E}[\theta \mid n_k]).
\end{aligned}
$$

The posterior expected value in this formula is found using Bayes' theorem. Given prior distribution function $P_0$ the posterior expected value is given by

$$
\mathrm{E}[\theta \mid n_k] = \frac{\int_0^1 \theta \, \mathrm{Bin}(n_k \mid k, \theta) \, \mathrm{d}P_0(\theta)}{\int_0^1 \mathrm{Bin}(n_k \mid k, \theta) \, \mathrm{d}P_0(\theta)}.
$$

In certain simple decision problems, we may wish to set our prediction method so as to maximize the probability of correct prediction (this would occur if the loss from our decision is purely dependent on whether our prediction is correct or not). To maximize the probability of correct prediction we use the prediction array which concentrates all of the predictive probability on those outcomes with the highest posterior expected value, conditional on the observed count. (There will always be a deterministic prediction method of this kind.)

Not all decision problems involving binomial prediction will involve maximizing the posterior probability of correct prediction. After all, strategies may arise based on other kinds of loss functions. To avoid the possibility of complete loss in a single gamble, Klotz (2000) discussed a betting strategy for multinomial prediction based on a loss function that penalises the use of 'extreme wagers' in which the entire bet is staked on a proper subset of the possible outcomes. We will not consider these kinds of alternative loss functions in the present analysis.

### 5. The frequent outcome approach

In many situations of interest leading to the binomial model, particularly those arising in gambling processes, the process generating the sequence $\boldsymbol{x}$ will have been designed to *attempt* to obtain a long-run proportion of $\theta = \frac{1}{2}$. If this can be achieved then the values in the observable series will be independent discrete uniform values (either zero or one) so that any prediction method is as good as any other.

Of course, the attempt to obtain a mechanism producing an equal long-run proportion of each outcome is often unsuccessful, and biases may exist within the process such that there are some disparities in these long-run proportions. If there is some nonzero probability of bias in the process but we have no *a priori* way to distinguish the likely direction of this bias, then this leads us to a judgment that the prior distribution will be nondegenerate and symmetric. This yields the following type of series.

**Definition 1.** If $\boldsymbol{x}$ is exchangeable and the distribution of $\theta$ is symmetric with $\mathrm{P}(\theta \neq \frac{1}{2}) > 0$ then we say that $\boldsymbol{x}$ is an *almost-uniform series*.

The almost-uniform series captures the notion of a series of values which are designed to be independent discrete uniform values, but which have some possibility of bias (i.e. some possibility of different long-run proportions of the outcomes). It is useful to note that the symmetry of the distribution of $\theta$ can be expressed equivalently as exchangeability of $(\theta, 1-\theta)$. This dual exchangeability assumption is sufficient to ensure that the series will either be a discrete uniform series (if the prior distribution is degenerate) or an *almost*-uniform series (if the prior distribution is nondegenerate).

O'Neill and Puza (2005) showed that these conditions for an almost-uniform series lead to a situation in which the probability of correct prediction is maximized by predicting whichever outcome has occurred *the most* in the observed data (if both outcomes have occurred an equal number of times then either prediction has an equivalent probability of success). This approach is known as the *frequent outcome approach*.

Suppose that we use the frequent outcome approach for prediction, with even odds when both outcomes occur the same number of times. (Hereafter, when we use the frequent outcome approach, we will mean that we are using this particular manifestation of this approach, with even odds used when both outcomes occur the same number of times. This avoids complicating our analysis later on.) In this case the elements in the prediction array are given by

$$\phi_{k,n_k} = \begin{cases} 0, & \text{if } n_k/k < \frac{1}{2}, \\ \frac{1}{2}, & \text{if } n_k/k = \frac{1}{2}, \\ 1, & \text{if } n_k/k > \frac{1}{2}. \end{cases}$$

With this prediction method we obtain

$$\pi_{k,n_k}(\theta) = \begin{cases} 1-\theta, & \text{if } n_k/k < \frac{1}{2}, \\ \frac{1}{2}, & \text{if } n_k/k = \frac{1}{2}, \\ \theta, & \text{if } n_k/k > \frac{1}{2}. \end{cases}$$

Taking the posterior expectation of the parameter $\theta$ we obtain the posterior probability of correct prediction. The specific form of this probability will depend on the form of the prior distribution. (It is common to use a beta prior distribution in the binomial model, since this is the conjugate of the binomial distribution.) However, if we are working with an almost-uniform series then this posterior probability of correct prediction using the frequent outcome approach, namely $\pi_{k,n_k}(\mathrm{E}[\theta \mid n_K])$, will be larger than the posterior probability of correct prediction for any prediction method not in accordance with this approach.

## 6. Accuracy function for specific processes

Given a specified prior distribution for $\theta$, the Bayesian method provides the means of determining the probability of correct prediction conditional on the observed values. Indeed, it is this probability of correct prediction that is used to determine the optimal prediction method under the Bayesian approach. This solves the prediction problem for a gambler who is able to make some prediction about the outcome of the next value in the observable series $x$.

While this information is all very useful for determining an optimal strategy, it gives us little understanding of the dynamics of applying this strategy to a random process with a specified level of bias (i.e. a specified value of the long-run proportion $\theta$). To assess this question we will instead need to look at the prior probability of correct prediction conditional on the parameter $\theta$.

In this section, we consider the probability of correct prediction conditional only on $\theta$. We use the *accuracy functions* $\pi_0, \pi_1, \pi_2, \ldots$ defined by

$$\pi_k(\theta) := \mathrm{P}(y_{k+1} = x_{k+1} \mid \theta) = \sum_{n=0}^{k} \pi_{k,n}(\theta) \operatorname{Bin}(n \mid k, \theta).$$

The accuracy function $\pi_k$ gives us the prior probability of correct prediction of $x_{k+1}$ assessed conditional on the parameter value $\theta$ (this is defined for any value $k$). Under the frequent outcome approach the elements in the summation for the accuracy function reduce down to

$$\pi_{k,n}(\theta) \operatorname{Bin}(n \mid k, \theta) = \begin{cases} \binom{k}{n} \theta^n (1-\theta)^{k-n+1}, & \text{if } n < k/2, \\ \frac{1}{2}\binom{k}{n} \theta^n (1-\theta)^{k-n}, & \text{if } n = k/2, \\ \binom{k}{n} \theta^{n+1}(1-\theta)^{k-n}, & \text{if } n > k/2. \end{cases}$$

The reader should be careful in interpreting what we are doing here. We are still proceeding on the basis that the prediction method is formed *with* knowledge of the previous observed values and *without* knowledge of the parameter. However, we then assess the probability of correct prediction under this prediction method *without* knowledge of the previous values and *with* knowledge of the parameter. This assessment would capture the situation of a gambler who sets his strategy using the observed values using the frequent outcome approach but then, *prior* to implementing it, asks himself how likely he is to correctly predict some future value *assuming* a particular level of bias in the process generating the values.

Since the accuracy function $\pi_k$ is a sum of polynomials of its argument value, it will be a polynomial function. We can obtain some useful expressions for this function by finding its recursive form, and then using this to express it in condensed and expanded polynomial form. We can also obtain a useful convergence result for the polynomial which lends additional credence to the frequent outcome approach. (Derivation of these results is in Appendix A.)

We let $C_0, C_1, C_2 \ldots$ be the Catalan numbers defined by $C_{i-1} = (2i - 2)!/i!\,(i - 1)!$ for all $i \in \mathbb{N}$ (see Koshy (2008)) and we define the functions $H_0, H_1, H_2, \ldots$ by

$$H_a(\theta) := \binom{2a}{a}\left(\frac{1}{2}\theta^a(1-\theta)^a - 2\theta^{a+1}(1-\theta)^{a+1}\right).$$

**Recursive form.** For all $a = 0, 1, 2, 3, \ldots$ we have the recursive relationship

$$\pi_{2a+2}(\theta) = \pi_{2a+1}(\theta) = \pi_{2a}(\theta) + H_a(\theta).$$

**Condensed polynomial form.** For all $a = 0, 1, 2, 3, \ldots$ we have

$$\pi_{2a+2}(\theta) = \pi_{2a+1}(\theta) = 1 - \sum_{i=1}^{a} \theta^i (1-\theta)^i C_{i-1} - 2\binom{2a}{a}\theta^{a+1}(1-\theta)^{a+1}.$$

**Expanded polynomial form.** For all $a = 0, 1, 2, 3, \ldots$ we have

$$\pi_{2a+2}(\theta) = \pi_{2a+1}(\theta) = 1 - \theta - \sum_{i=1}^{a+2} \alpha_{a,i}\theta^{a+i}.$$

(The coefficient values $\alpha_{a,i}$ are set out in Appendices A and B.)

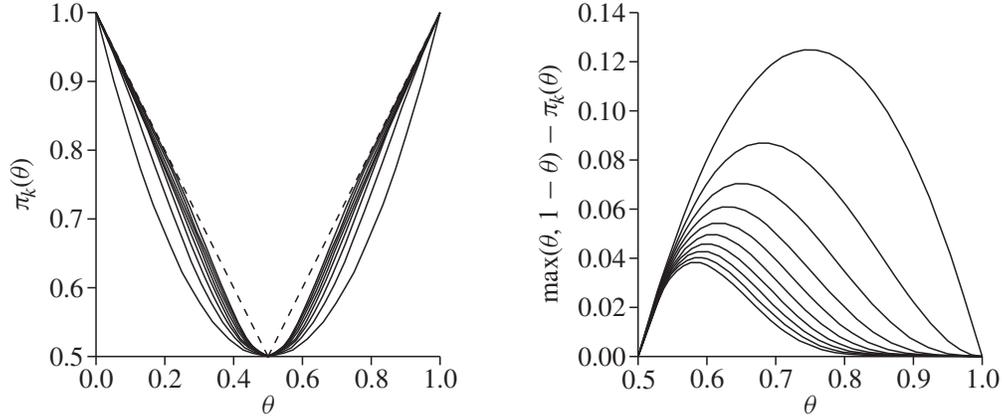

FIGURE 1: Accuracy functions converging to the idealised prediction probability.

**Convergence result.** We have $\lim_{k \to \infty} \pi_k(\theta) = \max(\theta, 1 - \theta)$.

The coefficients for the expanded polynomial form for the accuracy function are given in Appendix B. The first ten functions are

$$\pi_2(\theta) = \pi_1(\theta) = 1 - 2\theta + 2\theta^2,$$

$$\pi_4(\theta) = \pi_3(\theta) = 1 - \theta - 3\theta^2 + 8\theta^3 - 4\theta^4,$$

$$\pi_6(\theta) = \pi_5(\theta) = 1 - \theta - 10\theta^3 + 35\theta^4 - 36\theta^5 + 12\theta^6,$$

$$\pi_8(\theta) = \pi_7(\theta) = 1 - \theta - 35\theta^4 + 154\theta^5 - 238\theta^6 + 160\theta^7 - 40\theta^8,$$

$$\pi_{10}(\theta) = \pi_9(\theta) = 1 - \theta - 126\theta^5 + 672\theta^6 - 1380\theta^7 + 1395\theta^8 - 700\theta^9 + 140\theta^{10}.$$

We can see from the above convergence result that $\lim_{k \to \infty} \pi_k(\theta) = \max(\theta, 1 - \theta)$. This limiting value is the probability of correct prediction that would occur if we were able to use the ideal prediction method that accrues when our prior belief about $\theta$ is certain knowledge of its true value (i.e. we already know the true value $\theta$ when we formulate our prediction method). To illustrate this convergence we plot the accuracy functions $\pi_1, \ldots, \pi_{20}$ and their differences from ideal prediction in Figure 1.

We can see from these plots that the convergence occurs quite rapidly for values of $\theta$ that are in the extremes of its range, and occurs more slowly for values of $\theta$ that are near (but not at) the centre of its range. We can also see that the value of $\theta$ which gives the greatest difference in accuracy between ideal prediction and the frequent outcome approach moves inward to the middle of the range as $k$ increases.

This is all intuitively reasonable. If the value $\theta$ is in the extremes of its range then there is more of a loss in predictive power from our lack of knowledge of its value, but the frequent outcome approach will quickly detect and take advantage of the resulting imbalance in the outcomes. Conversely, if the value $\theta$ is near the middle of its range then there is less of a loss in predictive power from our lack of knowledge of its value, but the frequent outcome approach will take longer to detect and take advantage of the resulting imbalance in the outcomes.

### 7. Assessing the rate of convergence to ideal prediction

One of the reasons for analyzing the prior probability of correct prediction conditional on the parameter in the model is that it allows us to assess the rate at which prediction is likely to improve (assessed *a priori*) as the number of observed values increases. In particular, we are able to assess the rate of convergence of the accuracy function to the limiting situation of ideal prediction.

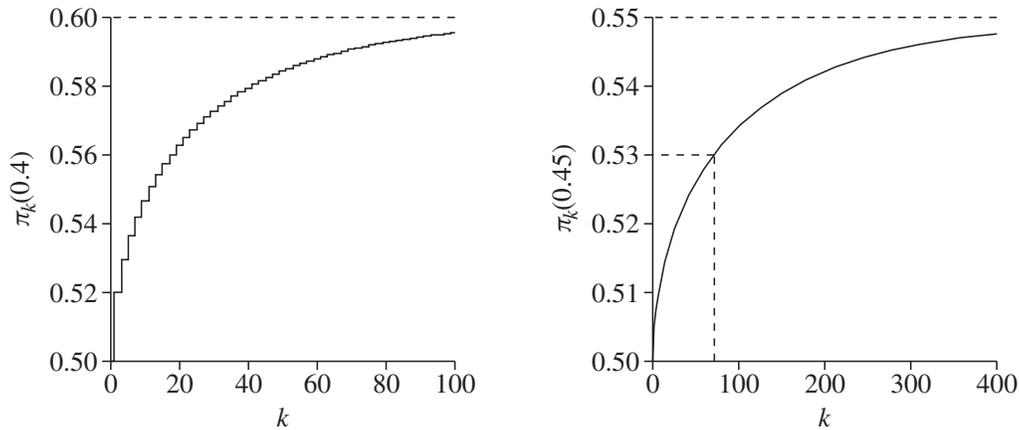

To do this, we take a fixed value $0 \leq \theta \leq 1$ and we observe $\pi_k(\theta)$ for $k = 1, 2, 3, \ldots$ to see convergence to its limiting value as $k \rightarrow \infty$. Plotting this as a step function we can see the rate at which the accuracy approaches the probability for ideal prediction.

Figure 2 illustrates the rate at which the accuracy of the frequent outcome approach converges to the accuracy of ideal prediction, where the parameter $\theta$ is known to be 0.4 and 0.45 respectively. Since this is assessed conditional on the parameter value it does not require specification of a prior distribution.

**Example 1.** A gambler wants to bet on the outcomes of a coin-tossing game in which the true long-run proportion of heads is 45% (i.e. it is a *very* biased coin). Since he judges the sequence of coin tosses to be an almost-uniform sequence, he decides to use the frequent outcome approach for his predictions. From the above convergence results we know that *in the limit* he will approach a 10% advantage over 'the house' with this strategy. That is, *in the long-run* he will win 55% of the time and the house will win 45% of the time. However, his actual advantage on any given prediction will always be less than a 10% advantage. How long, assessed *a priori*, will the gambler have to observe the outcomes in order to gain a 6% advantage in his prediction?

To get a 6% advantage over the house he needs a 53% chance of correct prediction. To see how long this takes we find the lowest value of $k$ for which $\pi_k(0.45) \geq 0.53$. Using the above formulae for the accuracy function we have

$$\pi_{70}(0.45) = 0.5298 \leq 0.53 \leq 0.5302 = \pi_{71}(0.45).$$

This means that, assessed *a priori*, it will take 71 coin tosses before we would expect there to be sufficient information for the gambler to gain a 6% advantage. This is marked in Figure 2.

One advantage of this kind of analysis is that it augments our knowledge of the optimal prediction method given by the frequent outcome approach. It gives us a sense of how long it takes for bias in a process designed to produce uniform outcomes to manifest itself with enough certainty for the frequent outcome approach to produce a substantial advantage in prediction. For gamblers dealing with processes of this kind (and for casinos trying to avoid advantages to gamblers), this kind of analysis will give them a sense of how long it takes for biases of a particular magnitude to give rise to advantages in prediction of a particular magnitude.

## 8. Extension to the multinomial model

The binomial model arises whenever we have a sequence of values that take on one of two possible values and are judged to be exchangeable. The natural extension to the multinomial model arises when the values in the sequence have some countable range which may be greater than two possible values.

For simplicity, the foregoing analysis has considered only the binomial model, and not its natural extension to the multinomial model. However, the frequent outcome approach is also the optimal method for multinomial prediction under an analogous set of assumptions regarding the prior distribution of the long-run proportions of the various outcomes (i.e. that the *vector* of long-run proportions is exchangeable).

Extending the accuracy function to accommodate multinomial prediction using the frequent outcome approach is conceptually simple, but gives rise to polynomial forms that are extremely cumbersome. Nevertheless, it is easy to extend the convergence proof in Appendix A.5 to demonstrate that the accuracy function for multinomial prediction also converges to the accuracy for ideal prediction given knowledge of the underlying parameter.

## 9. Concluding remarks

Previous research by O'Neill and Puza (2005) has already shown that the frequent outcome approach is the optimal prediction method for almost-uniform sequences (including the extension to the multinomial case) in the sense that it maximizes the posterior probability of correct prediction. Here we have shown that the accuracy of this method converges to idealized prediction given knowledge of the parameter value as the number of observations tends to infinity. This is a rather elegant and useful result. It means that if bias exists in the process under consideration, then the frequent outcome approach will find the most favourable outcome given enough observations and will eventually become effectively as good as a prediction method in which we already know the value of the parameter of interest.

The present results are useful in the context of prediction games involving almost-uniform sequences, such as those that commonly arise in simple gambling games. Our results show us how long a gambler should expect (*a priori*) to have to observe a biased process in order to raise his probability of correct prediction to some level that is within the bounds of idealized prediction.

## Appendix A. Derivation of equations for accuracy function

### A.1. The accuracy function

In this appendix we show the derivation of the various forms for the accuracy function in the main body of the paper. Recall from the body of the paper that

$$\pi_k(\theta) = \sum_{n=0}^{k} \pi_{k,n}(\theta) \operatorname{Bin}(n \mid k, \theta),$$

where

$$\pi_{k,n}(\theta) = \begin{cases} 1 - \theta, & \text{if } n < k/2, \\ \frac{1}{2}, & \text{if } n = k/2, \\ \theta, & \text{if } n > k/2. \end{cases}$$

We first use this formula to derive the recursive equation for the probability of correct prediction and then use this to derive the standard polynomial forms.

### A.2. Derivation of recursive form

For all $k \geq 0$ we define the functions $T_{k,n}$ by

$$T_{k,n}(\theta) := \begin{cases} \pi_{k,n}(\theta) \operatorname{Bin}(n \mid k, \theta), & \text{if } 0 \leq n \leq k, \\ 0, & \text{otherwise.} \end{cases}$$

To proceed, we define the functions $\overleftarrow{W}_{k,n}$ and $\overrightarrow{W}_{k,n}$ by

$$\overleftarrow{W}_{k,n}(\theta) := \begin{cases} 2(1-\theta)^2, & \text{if } n = k/2, \\ 1 - \theta, & \text{if } n \neq k/2, \end{cases} \qquad \overrightarrow{W}_{k,n}(\theta) := \begin{cases} 2\theta^2, & \text{if } n = k/2, \\ \theta, & \text{if } n \neq k/2. \end{cases}$$

From the defined values, it is easy to see that for all $k \geq 0$ we have

$$T_{k+1,n+1}(\theta) = \overrightarrow{W}_{k,n}(\theta) T_{k,n}(\theta) + \overleftarrow{W}_{k,n+1}(\theta) T_{k,n+1}(\theta).$$

From this relationship we obtain the recursive equation

$$\begin{aligned}
\pi_{2a+1}(\theta) &= \sum_{n=0}^{2a+1} T_{2a+1,n}(\theta) \\
&= \sum_{n=0}^{2a+1} (\overrightarrow{W}_{2a,n-1}(\theta) T_{2a,n-1}(\theta) + \overleftarrow{W}_{2a,n}(\theta) T_{2a,n}(\theta)) \\
&= \sum_{n=0}^{2a+1} (\theta T_{2a,n-1}(\theta) + (1-\theta) T_{2a,n}(\theta)) + (\theta(2\theta-1) + (1-\theta)(1-2\theta)) T_{2a,a}(\theta)
\end{aligned}$$

$$= \theta \sum_{n=0}^{2a} T_{2a,n}(\theta) + (1-\theta) \sum_{n=0}^{2a} T_{2a,n}(\theta) + (1-2\theta)^2 T_{2a,a}(\theta)$$

$$= \sum_{n=0}^{2a} T_{2a,n}(\theta) + (1-4\theta(1-\theta)) T_{2a,a}(\theta)$$

$$= \pi_{2a}(\theta) + \binom{2a}{a} \left( \frac{1}{2}\theta^a(1-\theta)^a - 2\theta^{a+1}(1-\theta)^{a+1} \right)$$

$$= \pi_{2a}(\theta) + H_a(\theta).$$

We also obtain the recursive equation

$$\pi_{2a+2}(\theta) = \sum_{n=0}^{2a+2} T_{2a+2,n}(\theta)$$

$$= \sum_{n=0}^{2a+2} (\overrightarrow{W}_{2a+1,n-1}(\theta) T_{2a+1,n-1}(\theta) + \overleftarrow{W}_{2a+1,n}(\theta) T_{2a+1,n}(\theta))$$

$$= \sum_{n=0}^{2a+1} (\theta T_{2a+1,n-1}(\theta) + (1-\theta) T_{2a+1,n}(\theta))$$

$$= \theta \sum_{n=0}^{2a+1} T_{2a+1,n}(\theta) + (1-\theta) \sum_{n=0}^{2a+1} T_{2a+1,n}(\theta)$$

$$= \sum_{n=0}^{2a+1} T_{2a+1,n}(\theta)$$

$$= \pi_{2a+1}(\theta),$$

which was to be shown.

## A.3.  Derivation of condensed polynomial form

Using repeated iteration of the above recursive equation, we have

$$\pi_{2a+2}(\theta) = \pi_{2a+1}(\theta)$$

$$= \pi_0(\theta) + \sum_{i=0}^{a} H_i(\theta)$$

$$= \pi_0(\theta) + \sum_{i=0}^{a} \binom{2i}{i} \left( \frac{1}{2}\theta^i(1-\theta)^i - 2\theta^{i+1}(1-\theta)^{i+1} \right)$$

$$= 1 - \sum_{i=1}^{a} \theta^i(1-\theta)^i \left( 2\binom{2i-2}{i-1} - \frac{1}{2}\binom{2i}{i} \right) - 2\binom{2a}{a}\theta^{a+1}(1-\theta)^{a+1}.$$

We can put this in terms of the Catalan numbers $C_0, C_1, C_2, \ldots$. To do this, we observe that

$$\binom{2a}{a} = \frac{2a!}{a!\,a!} = (a+1)\frac{2a!}{a!\,(a+1)} = (a+1)C_a$$

and

$$2\binom{2i-2}{i-1} - \frac{1}{2}\binom{2i}{i} = \frac{2i}{i}\frac{(2i-2)!}{(i-1)!\,(i-1)!} - \frac{i}{2i}\frac{(2i)!}{i!\,i!}$$

$$= \frac{(2i-2)!}{i!\,(i-1)!}2i - \frac{(2i-2)!}{i!\,(i-1)!}(2i-1)$$

$$= \frac{(2i-2)!}{i!\,(i-1)!}$$

$$= C_{i-1}.$$

From these equations we obtain

$$\pi_{2a+2}(\theta) = \pi_{2a+1}(\theta) = 1 - \sum_{i=1}^{a}\theta^i(1-\theta)^i C_{i-1} - 2(a+1)C_a\theta^{a+1}(1-\theta)^{a+1},$$

which was to be shown.

### A.4. Proof of convergence result

From the above, we have

$$\pi_{2a+2}(\theta) = \pi_{2a+1}(\theta) = 1 - \sum_{i=1}^{a}\theta^i(1-\theta)^i C_{i-1} - 2(a+1)C_a\theta^{a+1}(1-\theta)^{a+1}.$$

We want to find the limit of this polynomial function as $a \to \infty$. To do this, we will find the limits of the two main terms in the formula using the properties of the Catalan numbers. Letting $G_n(z) = \sum_{k=0}^{n} C_k z^k$, we have the generating function

$$G(z) := \lim_{n\to\infty} G_n(z) = \sum_{k=0}^{\infty} C_k z^k = \frac{2}{1+\sqrt{1-4z}}.$$

(For details of the properties of the Catalan numbers, including their generating function, see Koshy (2008).) From this equation we obtain the limit

$$\lim_{a\to\infty}\sum_{i=1}^{a}\theta^i(1-\theta)^i C_{i-1} = \theta(1-\theta)\lim_{a\to\infty}\sum_{i=0}^{a-1}\theta^i(1-\theta)^i C_i$$

$$= \theta(1-\theta)\lim_{a\to\infty} G_{a-1}(\theta(1-\theta))$$

$$= \theta(1-\theta)G(\theta(1-\theta))$$

$$= \frac{2\theta(1-\theta)}{1+\sqrt{1-4\theta(1-\theta)}}$$

$$= \frac{2\theta(1-\theta)}{1+\sqrt{(1-2\theta)^2}}$$

$$= \frac{2\theta(1-\theta)}{1+|1-2\theta|}$$

$$= \min(\theta, 1-\theta).$$

Now, using Sterling's approximation it is possible to show that $C_a \to 4^a/\sqrt{\pi a^3}$, in the sense that the ratio of these quantities tends to unity as $a \to \infty$. We therefore have

$$\lim_{a \to \infty} 2(a+1)C_a \theta^{a+1}(1-\theta)^{a+1}$$

$$= 2\theta(1-\theta) \lim_{a \to \infty} (a+1)C_a \theta^a (1-\theta)^a$$

$$= 2\theta(1-\theta) \lim_{a \to \infty} \frac{a+1}{a} C_a \frac{\sqrt{\pi a^3}}{4^a} \frac{(4\theta(1-\theta))^a}{\sqrt{\pi a}}$$

$$= 2\theta(1-\theta) \left( \lim_{a \to \infty} \frac{a+1}{a} \right) \left( \lim_{a \to \infty} C_a \Big/ \frac{4^a}{\sqrt{\pi a^3}} \right) \left( \lim_{a \to \infty} \frac{(4\theta(1-\theta))^a}{\sqrt{\pi a}} \right)$$

$$= 2\theta(1-\theta) \times 1 \times 1 \times 0$$

$$= 0.$$

(The last step here follows from the fact that $0 \le 4\theta(1-\theta) \le 1$ for all $0 \le \theta \le 1$.) Putting these two limits together we obtain

$$\lim_{k \to \infty} \pi_k(\theta) = 1 - \min(\theta, 1-\theta) = \max(\theta, 1-\theta),$$

which was to be shown.

### A.5. Alternative proof of convergence result

It is also quite simple to prove the convergence result without using the polynomial form for the accuracy function (relying solely on the law of large numbers). Since this is quite instructive, we also include this alternative proof.

If we use the frequent outcome approach for prediction then we have

$$\pi_k(\theta) = P(x_{k+1} = y_{k+1} \mid \theta)$$

$$= (1-\theta) P\left( \frac{n_k}{k} < \frac{1}{2} \;\Big|\; \theta \right) + \frac{1}{2} P\left( \frac{n_k}{k} = \frac{1}{2} \;\Big|\; \theta \right) + \theta P\left( \frac{n_k}{k} > \frac{1}{2} \;\Big|\; \theta \right).$$

For $\theta = \frac{1}{2}$ we have $\pi_k(\theta) = \frac{1}{2}$ so that the convergence result is trivial. We will therefore proceed for $\theta \ne \frac{1}{2}$. The weak law of large numbers means that $n_k/k$ converges in probability to $\theta$ as $k \to \infty$. This means that, for all $\varepsilon > 0$ we have

$$\lim_{k \to \infty} P\left( \theta - \varepsilon < \frac{n_k}{k} < \theta + \varepsilon \;\Big|\; \theta \right) = 1.$$

In particular, if $0 < \varepsilon < |\theta - \frac{1}{2}|$ then the above probability statement implies the following.

- If $\theta < \frac{1}{2}$ then

$$\lim_{k \to \infty} P\left( \frac{n_k}{k} < \frac{1}{2} \;\Big|\; \theta \right) = 1,$$

  so that $\lim_{k \to \infty} \pi_k(\theta) = 1 - \theta$.

- If $\theta > \frac{1}{2}$ then

$$\lim_{k \to \infty} P\left( \frac{n_k}{k} > \frac{1}{2} \;\Big|\; \theta \right) = 1,$$

  so that $\lim_{k \to \infty} \pi_k(\theta) = \theta$.

This means that $\lim_{k \to \infty} \pi_k(\theta) = \max(\theta, 1-\theta)$, which was to be shown.

## A.6. Expanded polynomial form

The condensed polynomial form included in the body of the paper is a simple way to represent the accuracy function. However, the expanded polynomial form is included here for completeness.

For all $a = 0, 1, 2, 3, \dots$ we have

$$\pi_{2a+2}(\theta) = \pi_{2a+1}(\theta) = 1 - \theta - \sum_{i=1}^{a+2} \alpha_{a,i}\theta^{a+i},$$

where the coefficients are given by

$$\alpha_{a,i} := 2(a+1)(-1)^{i-1}\binom{a+1}{t-1}C_a$$
$$+ \mathbf{1}_{\{i \le a\}} \sum_{j \le a} (-1)^{a+i-j}\binom{j}{a+i-j}C_{i-1} - \mathbf{1}_{\{a=0, i=1\}}.$$

## A.7. Derivation of expanded polynomial form

To derive the expanded polynomial form we will collect equivalent powers of $\theta$ in the condensed polynomial form. To do this we define

$$W_{i,j} := \begin{cases} (-1)^j \binom{i}{j} C_{i-1}, & \text{for } 0 \le j \le i, \\ 0, & \text{otherwise.} \end{cases}$$

Using the binomial theorem we then have the following useful equation:

$$C_{i-1}(1-\theta)^i = C_{i-1} \sum_{j=0}^{i} \binom{i}{j}(-1)^j \theta^j = \sum_{j=0}^{i} W_{i,j}\theta^j.$$

Using this equation we have

$$\pi_{2a+2}(\theta) = \pi_{2a+1}(\theta)$$
$$= 1 - \sum_{i=1}^{a} \theta^i (1-\theta)^i C_{i-1} - 2(a+1)C_a\theta^{a+1}(1-\theta)^{a+1}$$
$$= 1 - \sum_{i=1}^{a} \sum_{j=1}^{i} W_{i,j}\theta^{i+j} - 2(a+1)\sum_{j=0}^{a+1} W_{a+1,j}\theta^{i+j+1}$$
$$= 1 - \sum_{i=1}^{a} \sum_{t=i}^{2i} W_{i,t-i}\theta^t - 2(a+1)\sum_{t=a+1}^{2a+2} W_{a+1,t-a-1}\theta^t$$
$$= 1 - \sum_{t=1}^{2a} \left( \sum_{i \le \min(t,a)} W_{i,t-i} \right)\theta^t - 2(a+1)\sum_{t=a+1}^{2a+2} W_{a+1,t-a-1}\theta^t.$$

Now, it can be shown that $\sum_{i \leq t} W_{i,t-i} = \mathbf{1}_{\{t=1\}}$. But this means that

$$\pi_{2a+2}(\theta) = \pi_{2a+1}(\theta)$$

$$= 1 - \sum_{t=1}^{2a}\bigg(\sum_{i \leq \min(t,a)} W_{i,t-i}\bigg)\theta^t - 2(a+1)\sum_{t=a+1}^{2a+2} W_{a+1,t-a-1}\theta^t$$

$$= 1 - \sum_{t=1}^{a}\bigg(\sum_{i \leq t} W_{i,t-i}\bigg)\theta^t - \sum_{t=a+1}^{2a}\bigg(\sum_{i \leq a} W_{i,t-i}\bigg)\theta^t - \sum_{t=a+1}^{2a+2} 2(a+1)W_{a+1,t-a-1}\theta^t$$

$$= 1 - \theta\,\mathbf{1}_{\{a>0\}} - \sum_{t=a+1}^{2a}\bigg(\sum_{i \leq a} W_{i,t-i}\bigg)\theta^t - \sum_{t=a+1}^{2a+2} 2(a+1)W_{a+1,t-a-1}\theta^t$$

$$= 1 - \theta\,\mathbf{1}_{\{a>0\}} - \sum_{t=1}^{a}\bigg(\sum_{i \leq a} W_{i,a+t-i}\bigg)\theta^{a+t} - \sum_{t=1}^{a+2} 2(a+1)W_{a+1,t-1}\theta^{a+t}$$

$$= 1 - \theta\,\mathbf{1}_{\{a>0\}} - \sum_{t=1}^{a+2}\alpha_{a,t}\theta^{a+t},$$

which was to be shown.

## Appendix B. Coefficient values for the expanded polynomial form

From Appendix A we know that, for all $a = 0, 1, 2, 3, \ldots,$ the accuracy function can be expressed in expanded polynomial form as

$$\pi_{2a+2}(\theta) = \pi_{2a+1}(\theta) = 1 - \theta - \sum_{i=1}^{a+2}\alpha_{a,i}\theta^{a+i}.$$

TABLE 1: Coefficient values for expanded polynomial form.

| $\alpha_{a,i}$ | $i=1$ | $i=2$ | $i=3$ | $i=4$ | $i=5$ | $i=6$ |
|---|---|---|---|---|---|---|
| $a=0$ | 1 | $-2$ | | | | |
| $a=1$ | 3 | $-8$ | 4 | | | |
| $a=2$ | 10 | $-35$ | 36 | $-12$ | | |
| $a=3$ | 35 | $-154$ | 238 | $-160$ | 40 | |
| $a=4$ | 126 | $-672$ | 1 380 | $-1 395$ | 700 | $-140$ |
| $a=5$ | 426 | $-2 904$ | 7 425 | $-10 010$ | 7 546 | $-3 024$ |
| $a=6$ | 1 716 | $-12 441$ | 38 038 | $-64 064$ | 64 428 | $-38 766$ |
| $a=7$ | 6 435 | $-52 910$ | 188 188 | $-380 016$ | 477 750 | $-383 460$ |
| $a=8$ | 24 310 | $-223 652$ | 906 984 | $-2 134 860$ | 3 220 140 | $-3 231 360$ |
| $a=9$ | 92 378 | $-940 576$ | 4 282 980 | $-11 511 720$ | 20 252 100 | $-24 387 792$ |
| $a=10$ | 352 716 | $-3 938 662$ | 19 896 800 | $-60 116 760$ | 120 830 424 | $-169 744 575$ |

| $\alpha_{a,i}$ | $i=7$ | $i=8$ | $i=9$ | $i=10$ | $i=11$ | $i=12$ |
|---|---|---|---|---|---|---|
| $a=5$ | 504 | | | | | |
| $a=6$ | 12 936 | $-1 848$ | | | | |
| $a=7$ | 192 060 | $-54 912$ | 6 864 | | | |
| $a=8$ | 2 158 728 | $-926 211$ | 231 660 | $-25 740$ | | |
| $a=9$ | 20 369 349 | $-11 655 930$ | 4 374 370 | $-972 400$ | 97 240 | |
| $a=10$ | 170 143 974 | $-121 721 600$ | 60 920 860 | $-20 318 298$ | 4 064 632 | $-369 512$ |

The coefficient values in this formula can be calculated using the formula in Appendix A. This formula becomes computationally intensive if done manually, but it can be processed by a computer almost instantaneously. The coefficient values for the functions are provided in Table 1 for $a = 0, 1, 2, \ldots, 10$.